# Universal Meshes for the Simulation of Brittle Fracture and Moving Boundary Problems

Maurizio M. Chiaramonte, Evan S. Gawlik, Hardik Kabaria and Adrian J. Lew

**Abstract** Universal meshes have recently appeared in the literature as a computationally efficient and robust paradigm for the generation of conforming simplicial meshes for domains with evolving boundaries. The main idea behind a universal mesh is to immerse the moving boundary in a background mesh (the universal mesh), and to produce a mesh that conforms to the moving boundary at any given time by adjusting a few of elements of the background mesh. In this manuscript we present the application of universal meshes to the simulation of brittle fracturing. To this extent, we provide a high level description of a crack propagation algorithm and showcase its capabilities. Alongside universal meshes for the simulation of brittle fracture, we provide other examples for which universal meshes prove to be a powerful tool, namely fluid flow past moving obstacles. Lastly, we conclude the manuscript with some remarks on the current state of universal meshes and future directions.

*Dedicated to Michael Ortiz on the occasion of his 60$^{th}$ birthday.*

## 1 Introduction

Predicting and understanding the behavior of a propagating fracture has applications in a broad spectrum of disciplines. Perhaps the most renowned are applications in civil, mechanical, and aerospace engineering for the safe design of structural and mechanical components. More recently, a new wave of interest in understanding fracture propagation has risen due to the insurgence of hydraulic fracturing for the recovery of shale gas, as well as for engineering geothermal reservoirs. Alongside hydraulic fracturing, the practice of abyssal sequestration [1] for the disposal of radioactive waste also necessitates numerical tools capable of predicting the behav-

―――――――――――――――
Dept. of Mechanical Engineering and Institute for Computational and Mathematical Engineering Stanford University, Email:{`mchiaram,egawlik,hardik,lewa`}`@stanford.edu`

ior of fluid driven fractures. Beyond engineering, the modeling of fracturing finds relevance in geophysics, for example for the prediction of ice-sheet separation and its effect on global climate. Due to the pervasive nature of fracture mechanics in many disciplines, there is a need for a deeper understanding of fracture evolution accounting for the three-dimensionality of the fracturing process. These applications motivate the current efforts towards the creation of robust and computationally efficient numerical methods to approximate the solutions of such fracture evolution models. discuss, but is intended to shine some light on the pervasive nature of fracture mechanics amongst many disciplines. It is therefore apparent that there exists a need amongst many communities for robust and computationally efficient numerical methods for the prediction of fracture propagation.

From the numerical standpoint, one of the crucial challenges faced in this particular class of problems is the approximation of the evolving displacement discontinuity, which is the focus of the work presented here. Several approaches have been proposed in the literature to address this challenge. Albeit a comprehensive literature review is beyond the scope of this manuscript, a very broad classification of the predominant classes of methods capable of handling the evolution of a few cracks can be arguably categorized into basis-enriching methods or mesh-conforming methods. Additionally it is worthwhile mentioning numerical methods to approximate solutions of regularized theories of fracture. These theories, by assigning a finite width to the fracture, circumvent the need to explicitly track the crack geometry. Some examples are phase field methods [2], and Michael Ortiz's own contributions on eigenfracture [3] and eigenerosion [4, 5], to name a few. Also worth mentioning are methods for situations in which massive fragmentation appears, such as the seminal contributions by Michael Ortiz based on cohesive elements [6, 7, 8, 9, 10].

Basis-enriching methods, such as the Extended (XFEM) [11, 12] and Generalized (GFEM) [13, 14] finite element methods, endow the finite dimensional subspace with discontinuous functions. These methods circumvent the need to accommodate the evolving displacement discontinuity in the domain subdivision by implicitly representing it through the discontinuous basis functions. Numerical integration can be rather challenging, and, for problems such as hydraulic fracturing, when coupled governing equations need to be solved on the crack faces, these methods fail to provide a quality subdivision of the crack geometry. An example of the latter is illustrated in Fig. 1.

Alternatively, conforming methods envision generating a subdivision which accommodates exactly the evolving crack geometry. By ensuring that the mesh for the domain always conforms to the crack path, any displacement discontinuity along the crack is easily introduced. While the idea is simple it is nonetheless powerful. The robustness of this class of methods is limited by the generation of a quality conforming subdivision, a process which can be computationally demanding and prone to failure. Some examples of this approach are locally re-meshing methods as encountered in [15, 16, 17] as well as *r*-adaptive procedures as proposed in [18, 19, 20, 21]. Related to the latter are finite element spaces with embedded discontinuities [22, 23].

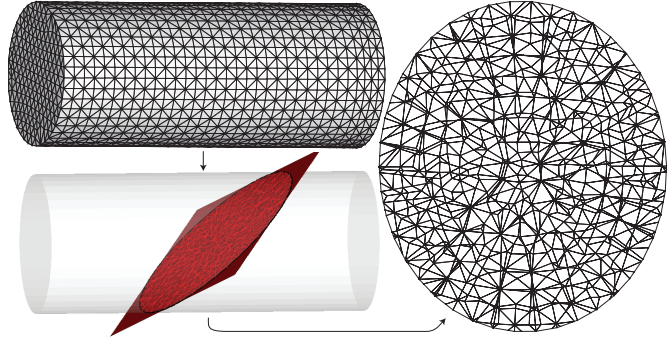

Fig. 1: An arbitrary cut (in red) through a quality tetrahedralization, representing an imaginary fracture, yields a poor discretization of the fracture faces.

Herein we present a few ideas for the simulation of brittle fracture that fall in the latter category of conforming mesh methods by taking advantage of *universal meshes*. Universal meshes is a paradigm for mesh generation that envisions the use of a single "background" mesh (the universal mesh) whose vertices closest to the crack geometry are perturbed to obtain a subdivision conforming to it. An example of such a perturbation is illustrated in Fig. 2. Because the same mesh can be deformed to conform to the geometry of a class of cracks, we say that the mesh is *universal* for such a class. The salient features of the method are its robustness, computational efficiency, and the mesh-independence of the solutions it provides (in fact, convergence).

We provide a description of the algorithm of universal meshes in § 2 followed by the presentation of the algorithm for the simulation of brittle fracture in § 3. Later, in § 4, we highlight some applications of universal meshes beyond brittle fracture. We conclude the manuscript on some recent developments of universal meshes in three dimensions in § 5.

## 2 Universal Meshes

We introduce the basic algorithmic ideas behind a universal mesh next. For concreteness, we focus the description on the aspects relevant to crack propagation, bearing in mind that similar ideas apply equally well to other classes of evolving domains, such as those encountered in fluid-structure interaction, as discussed later in § 4.

## *2.1 Algorithm*

To illustrate the discretization of an evolving domain with a universal mesh, we consider in this section the problem of triangulating a domain $\Omega(t) \subset \mathbb{R}^2$, $0 \leq t \leq T$, which contains an evolving crack. In other words,

$$\Omega(t) = \mathscr{D} \setminus \Gamma(t)$$

where $\mathscr{D} \subset \mathbb{R}^2$ is an open, bounded, polygonal domain and $\Gamma(t) \subset \mathscr{D}$ is a simple open rectifiable curve.

Let $\mathscr{T}_h$ be a triangulation of $\mathscr{D}$, hereafter referred to as the *universal mesh*. We use $h$ to denote the maximum diameter of an element of $\mathscr{T}_h$. We do not assume that the universal mesh conforms to $\Gamma(t)$ at any given time; in general, $\Gamma(t)$ may cut through elements of $\mathscr{T}_h$ arbitrarily, as in Fig. 2. Intuition would suggest, however, that a conforming mesh can be constructed by adjusting a few elements of $\mathscr{T}_h$ in a neighborhood of $\Gamma(t)$. This is the basic observation behind universal meshes.

To construct such a conforming mesh from $\mathscr{T}_h$, the following algorithm is adopted. First, a subset of edges in $\mathscr{T}_h$ lying near $\Gamma(t)$ is identified. We denote the union of these edges $\Gamma_h(t)$. Next, these edges are mapped onto $\Gamma(t)$ via the closest point projection $\pi : \mathbb{R}^2 \to \Gamma(t)$, with a suitable modification that places nodes precisely at the crack tips. Finally, the positions of nearby nodes are adjusted via a *relaxation* step that ensures the quality of the resulting triangulation.

The precise choices for the edges constituting $\Gamma_h(t)$ and the nodal adjustments adopted during relaxation are detailed in [24]. Briefly, $\Gamma_h(t)$ consists of *positive edges* of *positively cut* triangles in $\mathscr{T}_h$. To define these notions, one designates an orientation (positive or negative) for points in a neighborhood of $\Gamma(t)$. A triangle in $\mathscr{T}_h$ is called *positively cut* if it has two nodes on the positive side of $\Gamma(t)$ and one on the negative side. An edge is then called a *positive edge* if it belongs to a positively cut triangle and its endpoints both lie on the positive side of $\Gamma(t)$. A minor modification to $\Gamma_h(t)$ is made if a triangle in $\mathscr{T}_h$ has three nodes on $\Gamma_h(t)$; see [24] for details.

A key feature of the algorithm summarized above is its robustness. That is, the algorithm returns a valid mesh, for both the crack and the domain, with the quality of the elements bounded from below independently of the mesh size. To guarantee the above, three mild conditions need to be satisfied: (1) the background mesh is sufficiently refined in a neighborhood of $\Gamma(t)$, (2) all positively cut triangles in $\mathscr{T}_h$ are acute, and (3) the curve $\Gamma(t)$ is sufficiently smooth. This statement was proved for a domain with $C^2$ boundary (no cracks) in [25, 26]. The numerical examples strongly suggest that this should also be possible for domains with interfaces, such as cracks.

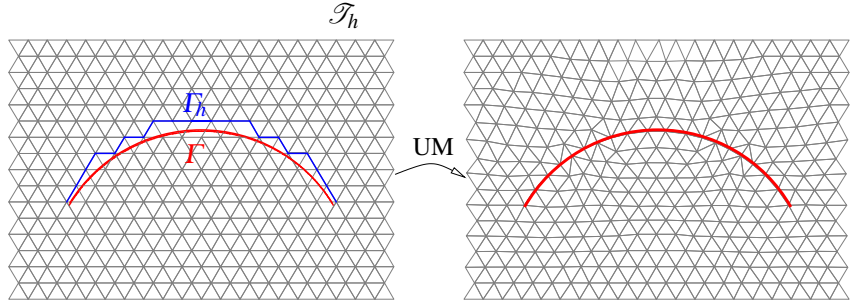

Fig. 2: Using a universal mesh, a triangulation conforming to a crack (right) is constructed by immersing the crack in a background triangulation $\mathcal{T}_h$ (left) and adjusting a few of its elements. This is accomplished by selecting a set of edges $\Gamma_h$ in the background triangulation that lie near the crack $\Gamma$, mapping them onto $\Gamma$ via the closest point projection, and relaxing a few nearby vertices to ensure the quality of the resulting triangulation.

## 3 Simulating Brittle Fracture with Universal Meshes

The obvious way in which a universal mesh is useful for the simulation of a propagating crack is by providing a mesh perfectly conforming to the crack at each step of its evolution. However, there are advantages that are less evident: the conforming mesh enables us to compute stress intensity factors to any order of accuracy, and the few mesh changes from step to step make it possible to retain much of the data structures in the computer implementation. The accuracy in the computation of the stress intensity factors is a determinant factor in observing convergence of the crack evolutions for "reasonable" mesh sizes.

In the following we present a numerical algorithm for the simulation of crack evolution in a restricted set of problems, as introduced in [24]. The presentation of the algorithm is followed by examples, including the formation of oscillatory crack paths in quenched plates, which requires very accurate stress intensity factors to converge.

### 3.1 Problem Statement

We consider the problem of *an always propagating* crack in an elastic medium, as defined next. We parametrize the crack evolution by the crack length $\ell \in [\ell_0, \ell_{max}]$ and we denote by $\mathscr{C}(\ell)$ the crack tip position for the crack of length $\ell$. Hence the crack set is given by $\mathscr{C}([\ell_0, \ell])$. The domain occupied by the cracked domain is denoted by $\Omega(\ell)$ and its boundary is decomposed into a portion over which dis-placements are prescribed, $\partial_d \Omega(\ell)$, and a portion over which boundary tractions are prescribed, $\partial_\tau \Omega(\ell)$. Further we let $\mathscr{C}([\ell_0, \ell]) \subseteq \partial_\tau \Omega(\ell)$.

The problem statement then reads: find the deformation $u(\cdot, \ell) : \Omega(\ell) \to \mathbb{R}^2$, the load scaling factor $C : [\ell_0, \ell_{max}] \to \mathbb{R}$, and the crack set $\mathscr{C}([\ell_0, \ell_{max}])$ such that the following holds for all $\ell \in (\ell_0, \ell_{\max}]$:

$$
\begin{aligned}
-\nabla \cdot (\mathbb{C} : \nabla u) &= b, & &\text{on } \Omega(\ell), \\
(\mathbb{C} : \nabla u) n &= t, & &\text{on } \partial_\tau \Omega(\ell), \\
u &= g, & &\text{on } \partial_d \Omega(\ell), \\
K_I[u] &= K_c, & & \\
K_{II}[u] &= 0, & & \\
\mathscr{C}([\ell_0, \ell^-]) &\subset \mathscr{C}([\ell_0, \ell]), & &\forall \ell^- < \ell,
\end{aligned}
$$

where $b(\ell) = C(\ell)\bar{b}(\ell), t(\ell) = C(\ell)\bar{t}(\ell), g(\ell) = C(\ell)\bar{g}(\ell)$. Here $\bar{b}(\cdot, \ell) : \Omega(\ell) \to \mathbb{R}^2$, $\bar{t}(\cdot, \ell) : \partial_\tau \Omega(\ell) \to \mathbb{R}^2$, and $\bar{g}(\cdot, \ell) : \partial_d \Omega(\ell) \to \mathbb{R}^2$ are arbitrary functions representing the "shape" of the body forces and boundary conditions. Effectively, for every crack length we know the "shape" of the applied body force ($b$), boundary tractions ($t$) and displacements ($g$), and we must solve for the linearly scaling coefficient $C(\ell)$ such that the condition $K_I[u] = K_c$ is always met, where $K_{I,II}[u]$ are the mode I and II stress intensity factors. The condition $K_{II}[u] = 0$ dictates the direction of crack propagation following the Principle of Local Symmetry [27].

The "always propagative crack" problem circumvents some of the more delicate issues in crack propagation, such as crack arrest and catastrophic crack propagation, regularity of the crack path, and competition among multiple cracks. The algorithm introduced next is applicable to this simpler class of problems.

## *3.2 Crack propagation algorithm*

There are three critical steps in the computation of the evolution of brittle crack paths: (1) the generation of a triangulation that conforms to the cracked domain, (2) the calculation of the elasticity fields, and (3) the evaluation of the stress intensity factors for curvilinear cracks. The steps are highlighted in Fig. 3.

We construct a triangulation that conforms to each cracked domains from a universal mesh, as described in §2. To ensure that the elasticity fields are sufficiently resolved, we draw on a class of finite element methods for domains with corners or cracks that retain optimal convergence rate for elements of any order in the face of singular solutions [28], in contrast to standard methods or methods with enrichments. Lastly, given that we compute with higher order solutions of the elasticity fields, we employ a family of *interaction integrals* specifically designed to handle curvilinear cracks [29] which yield stress intensity factors that converge rapidly to the exact ones (namely, they converge with twice the rate of convergence of the derivatives of the solution of the elasticity fields; a motivation to the use higher-order finite element methods).

We approximate the crack set $\mathscr{C}([\ell_0, \ell])$ by a cubic spline interpolant $\Gamma_\ell^h$ of a finite set of crack tips $\mathscr{A}_\ell = \{x_n\}_{n=0}^{(\ell-\ell_0)/\Delta\ell}$, where $\ell_0$ is the initial crack length and $\Delta\ell > 0$ is a crack discretization parameter. For a discrete crack $\Gamma_\ell^h$, $\ell$ indicates the chord length (the length along the polygonal line formed by points in $\mathscr{A}_\ell$, plus the initial crack length $\ell_0$) instead of its length. At any value of $\ell \geq \ell_0$, the algorithm proceeds through the following three steps:

1. Generate a conforming triangulation to the crack $\Gamma_\ell^h$.
2. Find $u^h(\ell) \approx u(\ell)$ and $C^h(\ell) \approx C(\ell)$.
3. Advance the crack in the direction $d(K_{II}^h[u^h]/K_I^h[u^h])$ by $\Delta\ell$ ( namely $\mathscr{A}_{\ell+\Delta\ell} = \mathscr{A}_\ell \cup \{x_{(\ell-\ell_0)/\Delta\ell} + d(K_{II}^h[u^h]/K_I^h[u^h])\Delta\ell\}$ ).

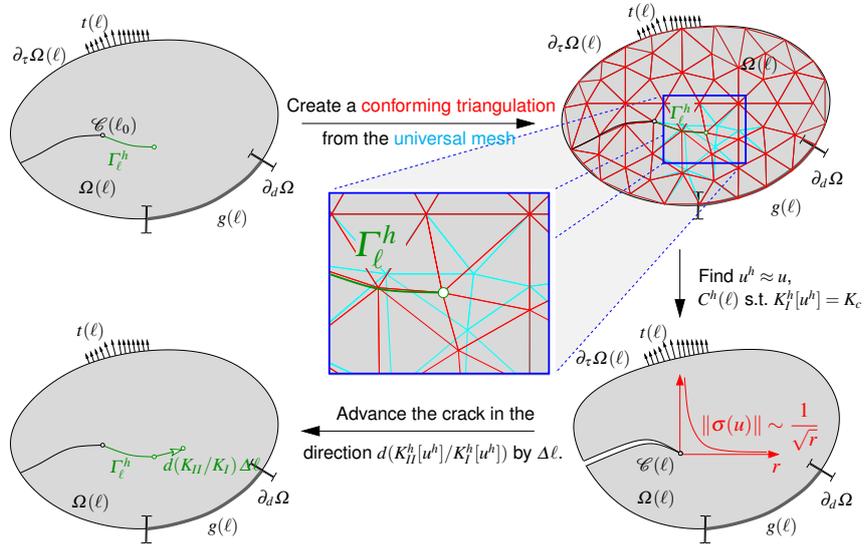

Fig. 3: The critical steps in the crack advancement algorithm.

Here the direction $d : \mathbb{R} \to S^1$ is chosen such that an infinitesimally short kink at the chosen angle satisfies the Principle of Local Symmetry ($K_{II} = 0$) up to first order in the kink angle itself. Evaluating $d$ in this way sidesteps the potentially computationally intensive alternative of explicitly solving for the direction $d$, but it likely restricts the order of convergence of the algorithm.

## *3.3 Examples*

We next showcase the application of the algorithm for the simulation of brittle fracture to two examples. The first example is a crack propagating along a circular arc, which we compare against an exact solution, and the second is a crack propagating in a perforated plate undergoing three-point bending, which we compare against experimental results. We also show (preliminary) results on the application of a slight modification of this algorithm to the computation of crack paths in a rapidly quenched plate [30].

### 3.3.1 A crack propagating along a circular arc

The displacement and stress fields of an infinite medium that contains a crack shaped as a circular arc subjected to far-field stresses and traction-free faces was computed in [31]. The corresponding stress intensity factors can be found in [32]. We use this solution to construct a loading history that, when is applied as Dirichlet boundary conditions to a square-shaped domain, as illustrated in Fig. 4a, causes the crack to propagate along a circular arc. For details on the construction of such a loading history we refer the interested reader to [24].

Fig. 4a shows a square-shaped domain $\Omega$ with a pre-existing crack of radius $R = 2$ and angular span $\vartheta_0 = \pi/8$. The analyses were carried out with four progressively refined universal meshes. The coarsest background mesh as well as the conformed mesh are shown in Fig. 5, and their refinements were constructed by recursively subdividing each triangle of the background mesh into four similar ones. The ratio $\Delta \ell / h \approx 2$ was kept constant over all simulations, where as usual $h$ denotes the maximum diameter of an element in a triangulation. As shown in Fig. 4b, the crack path converges to a circular arc, and the convergence curves for the $L^p([\ell_0, \ell_{max}])$, $p = 2, \infty$ and $H^1([\ell_0, \ell_{max}])$ norms are shown in Fig. 4c. Notably, convergence of the tangents to the crack path is also obtained.

This simple, nonetheless illustrative example, suggests that the algorithm is indeed convergent, and hence that the computed paths are largely independent of the chosen mesh.

### 3.3.2 Perforated plate

We next present the problem of a perforated plate undergoing three-point-bending. The problem setup is illustrated in Fig. 6a. We performed the simulations for three configurations of the initial crack position ($d$) and length ($\ell_0$). The values are tabulated in Fig. 6a. In Fig 6b we illustrate a universal mesh employed for one of the three simulations. It is worthwhile to note the adaptive nature of the background triangulation; in fact universal meshes can be easily adopted in conjunction with adaptive refinement. For each of the three simulations we generated a slightly different universal mesh to comply with the varying location of the initial crack.

Universal Meshes for Brittle Fracture and Moving Boundary Problems      9

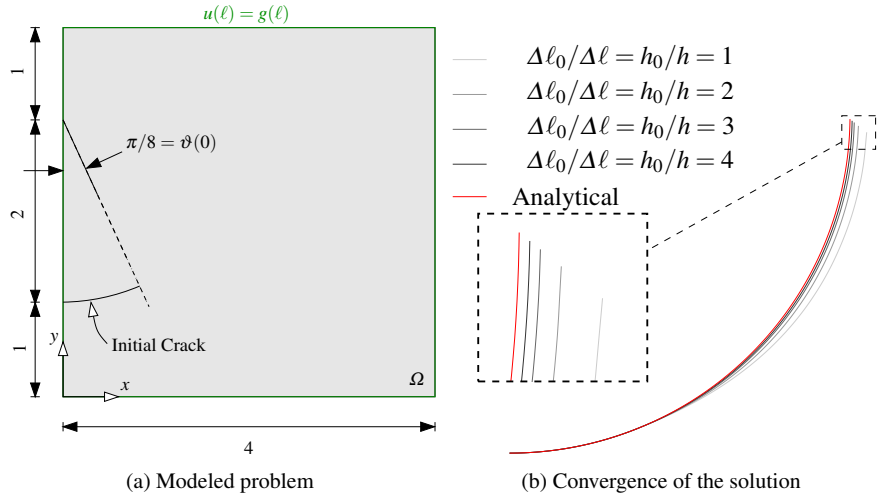

(a) Modeled problem      (b) Convergence of the solution

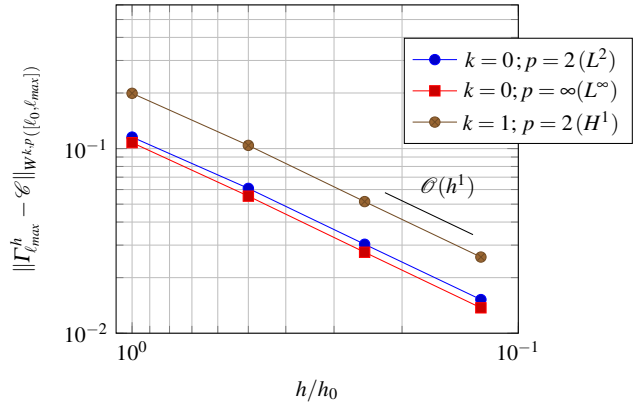

(c) Convergence of the crack path

Fig. 4: The circular arc crack problem.

Experimental results for this test setup are available in [33, 15]. The experiments were performed on polymethyl methacrylate (PMMA) plates. A comparison of the computed crack paths with digitized points from [33, 15] show a good agreement with experimental results. Further, relative convergence studies were performed on the computed crack paths, and the results, that can be found in [24], show a similar behavior to the one observed in Fig. 4.

10 Maurizio M. Chiaramonte, Evan S. Gawlik, Hardik Kabaria and Adrian J. Lew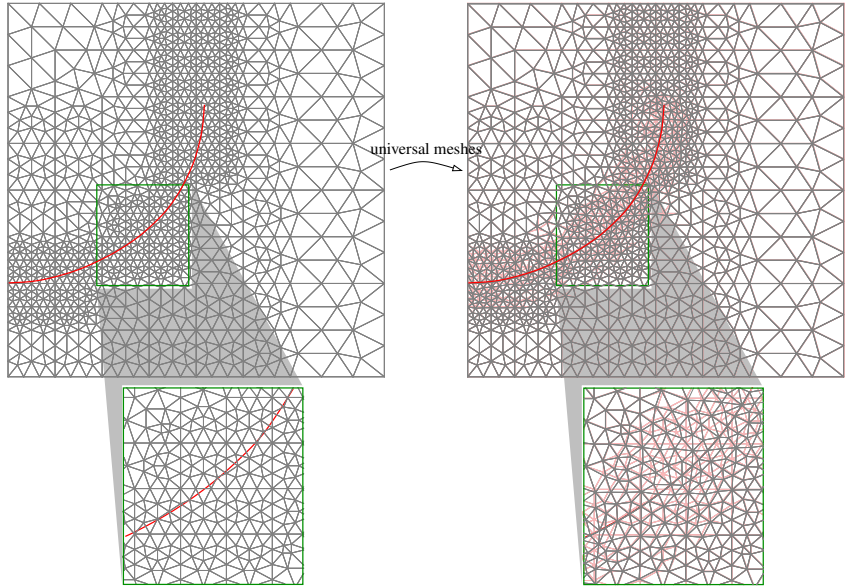

Fig. 5: The universal mesh used for the entire simulation (left) and the conformed triangulation at the end of the simulation with the background mesh in red (right).

### 3.3.3 Crack path instabilities in a quenched plate

Lastly we concisely present the problem of a thermally driven fracture in a quenched plate. The problem consists of a plate of finite width cracked along its center line, with low toughness ($K_c$), that, after being heated to temperature $\theta^+$, is immersed in an ice bath at temperature $\theta^-$ with a constant velocity ($v$). Refer to Fig. 8 to supplement the above description.

Depending on the material parameters, the presence of small deviations from the idealized descriptions above, and the configuration of the experiment, the crack path is expected to develop oscillations. Fig. 9 showcases the results of experiments performed by Yuse and Sano [34]. In Fig 10 we showcase some snapshots of the computed crack path for one set of inputs. These were computed through a modification to time-dependent problems of the algorithm introduced here. Details will appear in [30].

Although not shown here, the crack paths are converged up to a small tolerance. In our experience, this problem benefitted immensely from the high-order computation of the stress intensity factors; our previous attempts with low-order methods required inordinately large meshes to begin displaying some form of mesh-independent results. We hope to use this platform to better understand the underlying physics.

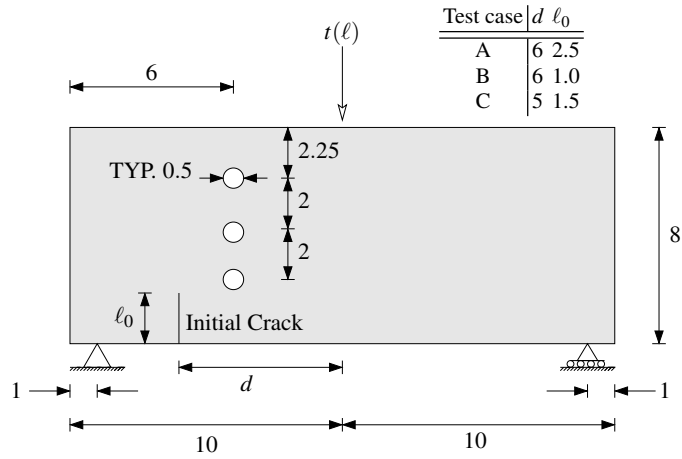

(a) Modeled domain

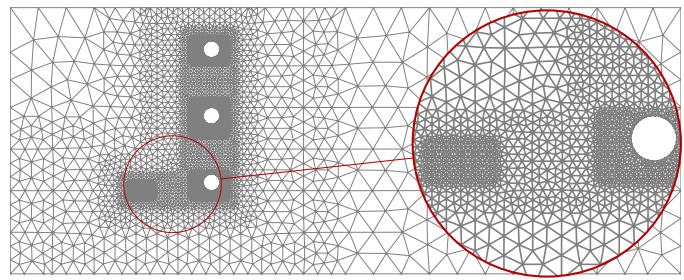

(b) Universal mesh

Fig. 6: Geometry for a plate with holes and the universal mesh adopted.

## 4 Beyond Brittle Fracture: Moving Boundary Problems

In addition to crack propagation, a variety of problems in science and engineering involve partial differential equations posed on domains that change with time. Such problems, collectively referred to as *moving-boundary problems*, appear in studies of fluid-structure interaction, phase-transitions, free-surface flows, aeroelasticity, and biolocomotion, to name a few. In this section, we demonstrate the applicability of universal meshes to this broader class of problems.

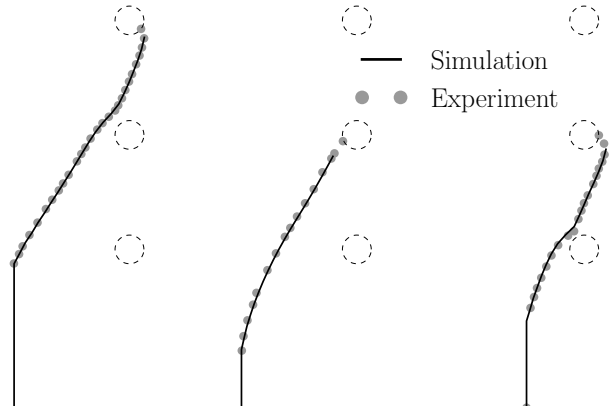

Fig. 7: Comparison between experimental results and computed crack paths. Experimental results were digitized from [33, 15].

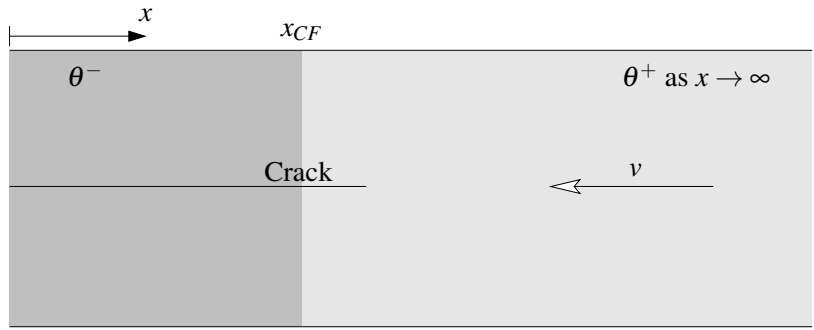

Fig. 8: Quenched plate problem setup.

## 4.1 Examples: Flow past moving obstacles

In the setting of fluid-structure interaction, universal meshes provide a conforming discretization of the fluid domain at all times, allowing finite element spaces of any desired order of accuracy to be used to spatially discretize the Navier-Stokes equations. This conforming discretization can be made to deform smoothly over time intervals that are short in comparison to the mesh spacing, thereby allowing standard numerical integrators to be used to solve the resulting system of ordinary differential equations. A projector (such as the nodal interpolation operator) is then used to transfer information between finite element spaces each time nodal positions change discontinuously. Details of this procedure are given in [44, 45], and rigorous theoretical bounds for a wide class of linear problems guarantee the high-order nature of the resulting numerical scheme when high-order finite elements are adopted [46, 47].

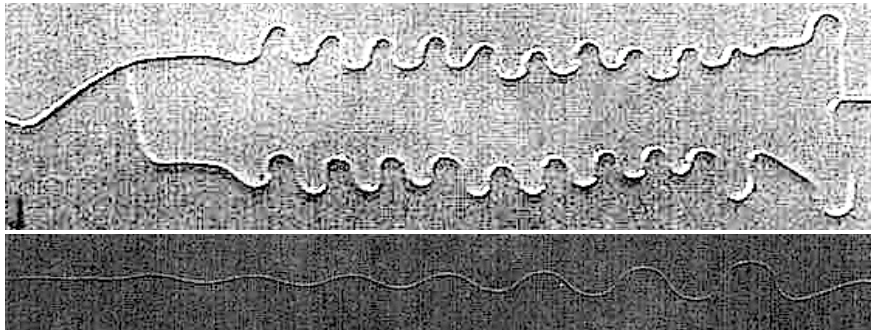

Fig. 9: Representative results of experiments of wavy crack patterns in rapidly quenched plates [34]. In both cases shown above, crack propagation along a straight crack is unstable. These cases correspond to different immersion speeds.

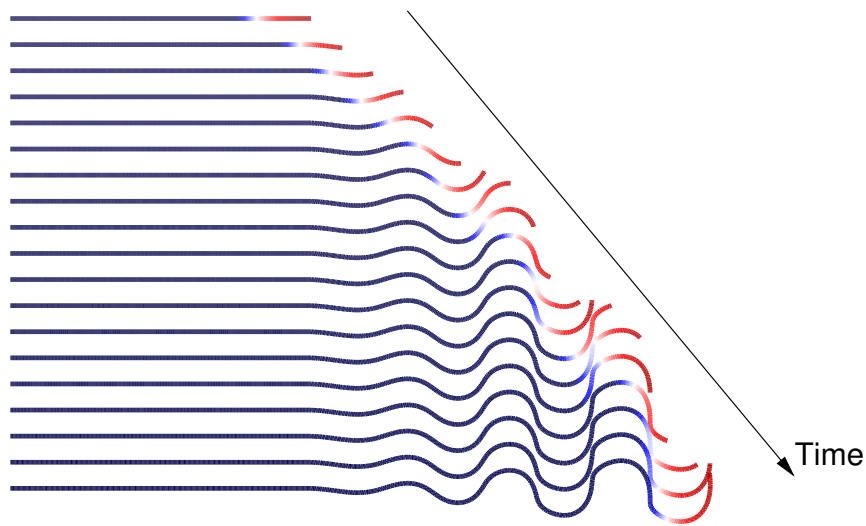

Fig. 10: Computed evolution of a thermally driven crack in a quenched plate. The contours show the temperature profile along the crack with dark blue representing $\theta^-$ and dark red representing $\theta^+$

As an example, we consider in Fig. 11 the solution of incompressible, viscous flow past a rotating propeller at Reynolds number $Re = 290$. We solved the problem using a universal mesh having adaptive refinement in a neighborhood of the propeller, together with Taylor-Hood finite elements. Figure 11 shows contours of the vorticity at two instants in time. The robust nature of the method is patent in this example, as traditional deforming-mesh methods could easily encounter difficulties with mesh entanglement upon rotation of the propeller.

As a second example, we consider in Fig. 11 the solution of incompressible, viscous flow past a pair of NACA0015 airfoils that change their pitch sinusoidally in time. We solved the problem using a universal mesh together with Taylor-Hood finite elements. For simplicity, the tips of the airfoils were blunted so that the algorithm described in Section 2 (which applies to smooth geometries) could be applied in its most basic form. Figure 11 shows contours of the vorticity at two instants in time corresponding to the maximum and minimum pitch ($17°$ and $-17°$, respectively) of the airfoils.

Finally, we consider the solution of incompressible, viscous flow past an oscillating disk with unit diameter at Reynolds number $Re = 185$. We solved the problem using a universal mesh having adaptive refinement near the disk (see Fig. 12a), together with Taylor-Hood finite elements [44]. The disk's motion was prescribed using a sinusoidally varying vertical displacement with amplitude 0.2 and frequency equal to 0.8 times the natural shedding frequency of a fixed disk of the same diameter. Fig. 12b shows a snapshot of the contours of the vorticity. Fig. 12c shows the observed convergence of the drag and lift coefficient time series under refinement of the mesh, which were computed via direct integration over the boundary of the disk.

## 5 Outlook

Clearly for a universal mesh to be useful in engineering practice, it needs to be able to handle evolving geometries in three-dimensions. We show next some incipient results in this direction.

### *5.1 Universal Meshes for Smooth Three-Dimensional Domains.*

The construction of a universal mesh in three dimensions follows the steps described in §2. Namely, given a smooth closed surface $\Gamma \subset \mathscr{D} \subset \mathbb{R}^3$ immersed in a mesh of tetrahedra $\mathscr{T}_h$, we first identify a set of faces $\Gamma_h$ in $\mathscr{T}_h$ that lie near $\Gamma$. These faces are then mapped onto $\Gamma$ via the closest point projection, and nearby nodes are adjusted via a relaxation step that ensures the quality of the resulting mesh.

In analogy with the algorithm presented in §2, $\Gamma_h$ is chosen as the union of *positive faces* of positively cut tetrahedra in $\mathscr{T}_h$. A tetrahedron in $\mathscr{T}_h$ is called *positively cut* if it has three nodes on the non-negative side of $\Gamma$ and one on the negative side. A face is then called a *positive face* if it belongs to a positively cut tetrahedron and all three of its vertices lie on the non-negative side of $\Gamma$. The closest point projection defines a one-to-one mapping between a positive face and its image on $\Gamma$ provided that the mesh size is small compared to the local radius of curvature, and more importantly, provided some dihedral angles in the mesh are acute [48].

Some examples that illustrate the use of a universal mesh in three dimensions are given in Figs. 13-14. In Fig. 13, two meshes of tetrahedra conforming to an elephant undergoing changes in its posture were obtained from a single universal mesh. In Fig. 14, the same procedure was used to construct conforming meshes of tetrahedra of a human upper airway.

### *5.2 Universal Meshes for Evolving Curves on Surfaces*

With an eye towards evolving crack fronts in three dimensions, we show next some early results on how a universal mesh can conform to a smooth curve drawn over a smooth surface, triangulating the interior of the curve over the surface, and conforming the tetrahedra to the surface and the mesh as well, see Fig. 15. To do so, a planar parametrization of the surface was constructed, and a smooth approximation of the given curve immersed in it. A conforming surface triangulation to the curve was then achieved by using a modification of the algorithm in two dimensions, not described here, and mapping the resulting planar triangulation back to the surface.

**Acknowledgements** This work was supported by the Office of Technology Licensing Stanford Graduate Fellowship to Maurizio M. Chiaramonte, the National Science Foundation Graduate Research Fellowship to Evan S. Gawlik, and the Franklin P. Johnson Jr. Stanford Graduate Fellowship to Hardik Kabaria. Adrian J. Lew acknowledges the support of National Science Foundation; contract/grant number CMMI-1301396.

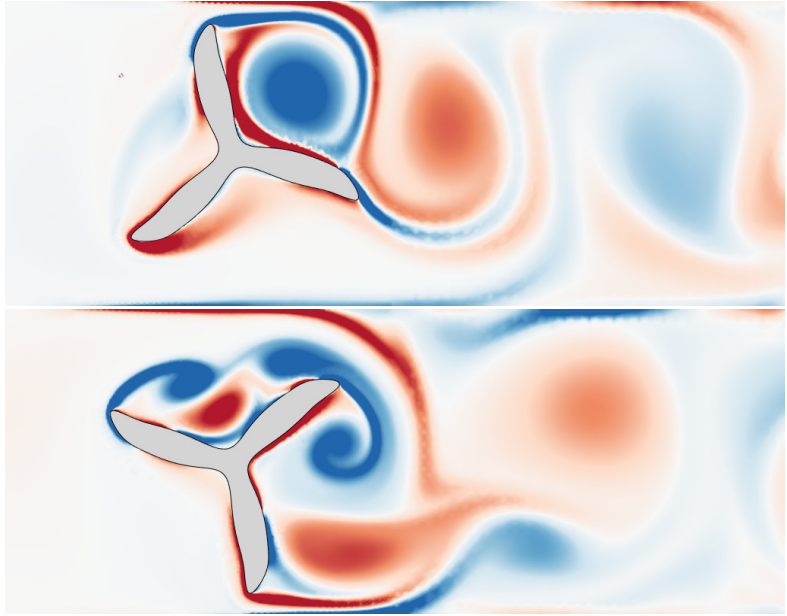

(a) Flow past a rotating propeller.

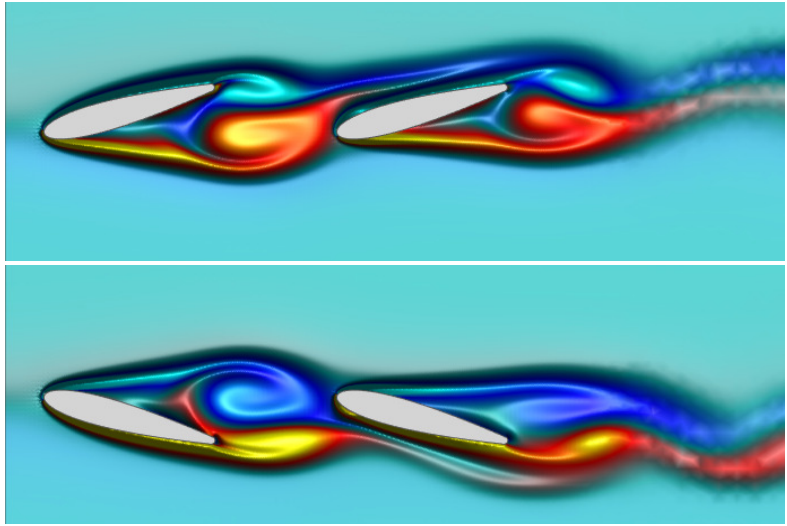

(b) Flow past pitching airfoils

Fig. 11: Vorticity contours for two representative examples of flow past a moving obstacle. The simulations consist of incompressible viscous flow, computed using a universal mesh together with Taylor-Hood finite elements.

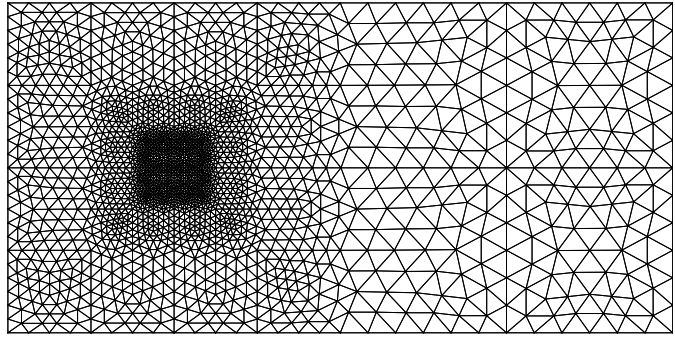

(a) Background mesh

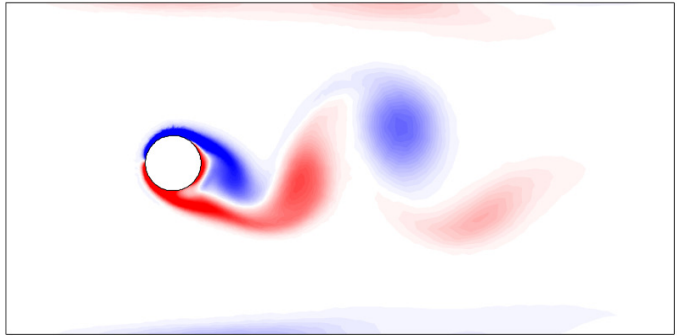

(b) Snapshot of the vorticity contours

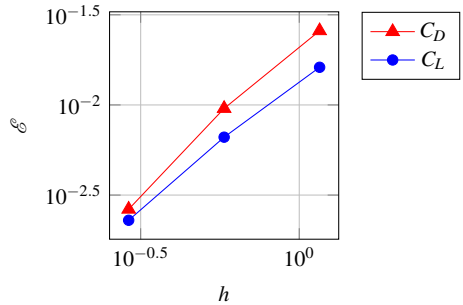

(c) Convergence of drag and lift coefficients

Fig. 12: Numerical simulation of incompressible, viscous flow past an oscillating disk using a universal mesh. In (a), the background mesh adopted during the simulation is shown. In (b), a snapshot of vorticity contours are shown. In (c), the convergence of the drag and lift coefficient time series under mesh refinement is shown. The reported error $\mathscr{E}$ is the square root of the integrated squared error $(C_i(t) - \bar{C}_i(t))^2$, $i = L, D$, over the time interval $[0, 1]$, relative to a reference solution $\bar{C}_i(t)$ obtained from a fine mesh with $h = 0.145$. Nearly quadratic convergence is observed. See [44] for details.

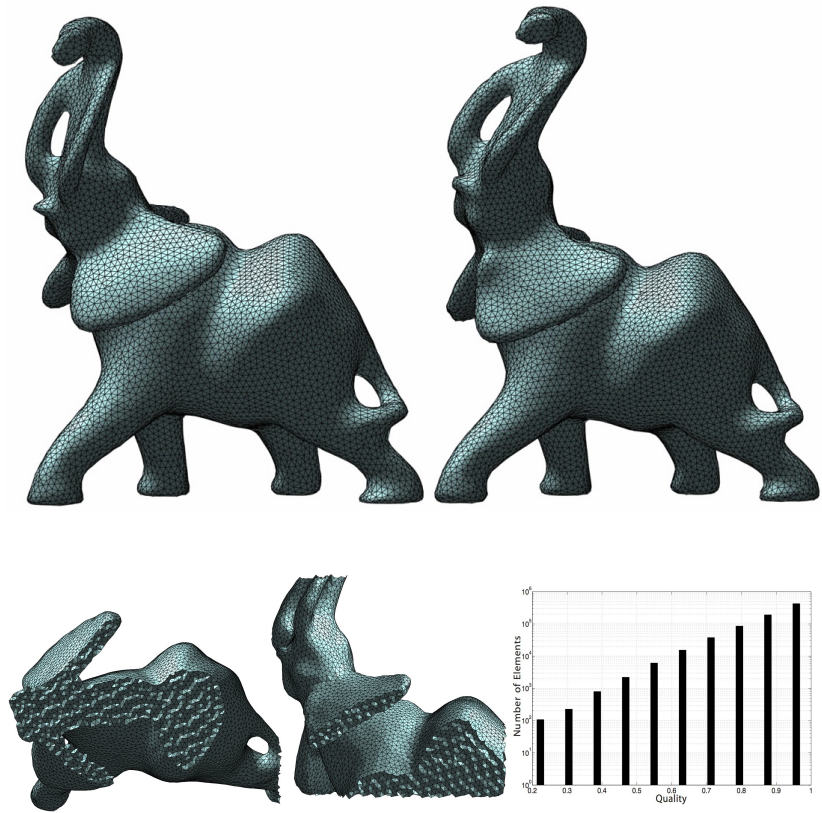

Fig. 13: The top row shows the tetrahedral meshes conforming to the geometry of an elephant in two different postures. They were obtained from the same background universal mesh, discarding exterior elements. In the bottom row, two cuts of the mesh displayed above and right are shown, together with the distribution of the quality of the tetrahedra in the same mesh, on a logarithmic scale.

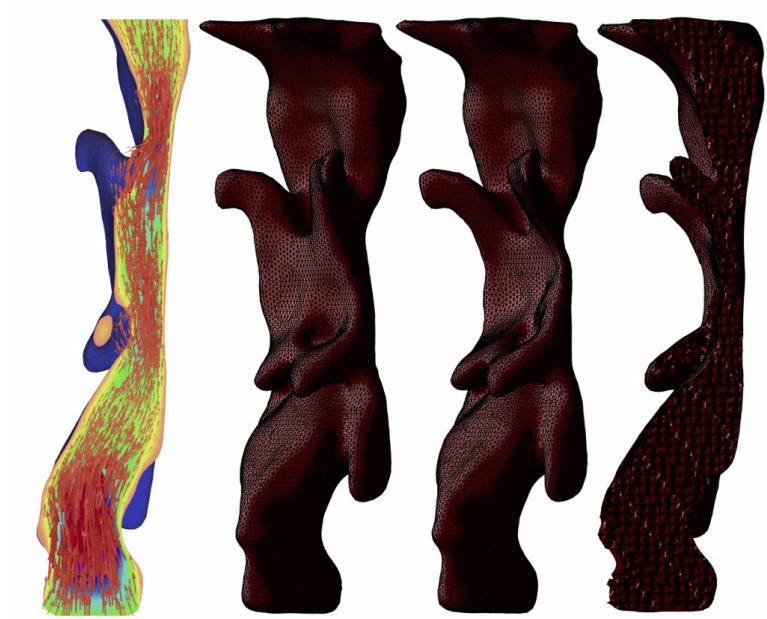

Fig. 14: The two figures in the middle show conforming meshes for two different configurations of a human upper airway. The figures on the far left and the far right show a cut through each tetrahedral mesh. We used a CT scan of a patient as an input for this example. The figure on the far left shows the contours of the velocity field computed by solving the Navier-Stokes equations inside. Coupling this tool to a solid mechanics analysis code for the upper airway would be useful to study the collapse of the upper airway, quite often the area of interest for patients diagnosed with obstructive sleep apnea.

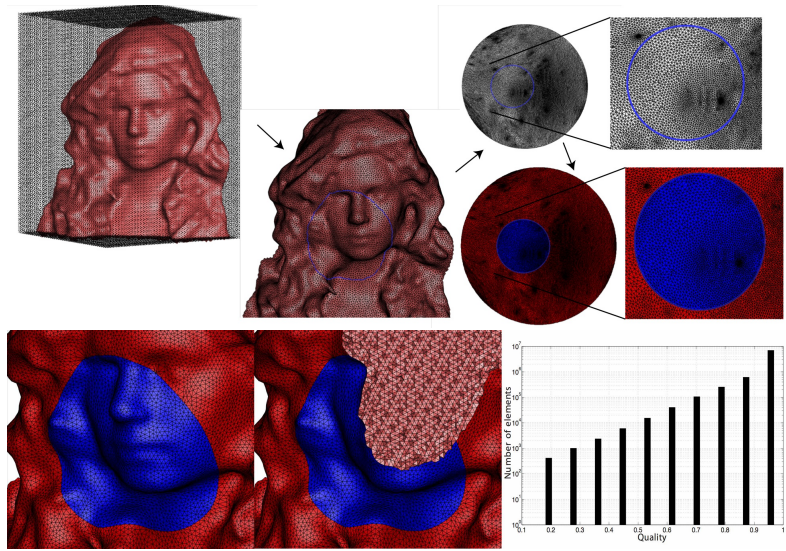

Fig. 15: A universal mesh can be used to generate tetrahedral meshes that conform to curves on surfaces. As a first step, a conforming mesh of the given surface is generated from the background universal mesh, as shown in the top left and top middle figures. We then construct a planar parametrization from the corresponding surface triangulation (top right), and conform the mapped mesh to the curve in the parametric planar domain. Mapping then back to the real space we obtain a surface triangulation that conforms to the curve on the surface (bottom left), and after a relaxation step of the nodes near the surface to ensure good quality of the tetrahedra, we obtain the resulting mesh (bottom middle). The distribution of the qualities of tetrahedra in this mesh is shown at the bottom right.